\documentclass[english]{article}
\usepackage{blindtext}
\usepackage[a4paper, total={6in, 8.5in}]{geometry}
\usepackage[T1]{fontenc}
\usepackage[latin9]{inputenc}
\usepackage{babel}
\usepackage{amsfonts,color}
\usepackage{graphicx}
\usepackage[font=small,labelfont=bf]{caption}

\newtheorem{dfn}{Definition}[section]
\newtheorem{theorem}[dfn]{Theorem}

\newtheorem{remark}[dfn]{Remark}

\usepackage{amssymb} 
\usepackage{amsmath}

 \global\long\def\sbr#1{\left[ #1\right] }
 \global\long\def\cbr#1{\left\{  #1\right\}  }
 
 \global\long\def\rbr#1{\left(#1\right)}

 \global\long\def\R{\mathbb{R}}
  \global\long\def\N{\mathbb{N}}
  \global\long\def\Q{\mathbb{Q}}
  \global\long\def\Z{\mathbb{Z}}

\global\long\def\M{\mathbb{M}}

 \global\long\def\dd#1{\textnormal{d}#1}

 \global\long\def\TTV#1#2#3{\text{TV}^{#3}\!\rbr{#1,#2}}

 \global\long\def\ra{\rightarrow}
 
 \global\long\def\ns{\infty}
 \global\long\def\8{\infty}

 \global\long\def\ns{\infty}

\global\long\def\cross#1#2#3{\text{n}^{#1}\!\rbr{#2,#3}}

\global\long\def\Ucrossemph#1#2#3{\emph{u}^{#1}\!\rbr{#2,#3}}
 \global\long\def\Dcrossemph#1#2#3{\emph{d}^{#1}\!\rbr{#2,#3}}
 \global\long\def\crossemph#1#2#3{\emph{n}^{#1}\!\rbr{#2,#3}}

\title{Quadratic variation and local times of the horizontal component of the Peano curve (square filling curve)}

\author{ Phumlani L. Zondi, Darlington Hove, Rafa{\l } M. {\L }ochowski\footnote{e-mail: rlocho@sgh.waw.pl, ORCID no.: 0000-0001-5427-8508 }, Farai J. Mhlanga \footnote{ORCID no.: 0000-0002-8775-9073} }

\begin{document}

\maketitle

\begin{abstract}
We show that the horizontal component of the Peano curve has quadratic
variation equal the limit of quadratic variations along the Lebesgue
partitions for grids of the form $3^{-n}p\Z+3^{-n}r$, $n=1,2,\ldots$, where $p$ is a rational number, while $r$ is irrational number, but the value of such quadratic variation depends on $p$.

This also yields that the horizontal component of the Peano curve is  an example of a deterministic function possessing local time (density of the occupation measure) with respect to the Lebesgue measure, whose
local time can be expressed as the limit of normalized numbers
of interval crossings by this function but the normalization is not a smooth function of the width of the intervals.

These two features distinct the horizontal component of the Peano curve from the trajectories of the Wiener process, which is widely used in financial models.
\end{abstract}

\section{Introduction}

The aim of this note is to investigate existence of quadratic
variation and local time of the $x$-component of the Peano curve.

The Peano curve is a well known continuous function ${\cal P}:[0,1]\ra[0,1]^{2}$,
which attains (as its value) any point of the unit square $[0,1]^{2}$,
that is, ${\cal P}$ is continuous and the image of the segment $[0,1]$ under
${\cal P}$ is the whole unit square $[0,1]^{2}$. Both of its components
$(x,y)$ (${\cal P}(t)=(x(t),y(t)),\ t\in[0,1]$) reveal fractal
structure, stemming from self-similarity. The formal definition of
the Peano curve is the following (see \cite[Chapt. 3]{Sagan:1994}). Let $t=0_{3}t_{1}t_{2}t_{3}\ldots:=\sum_{j=1}^{+\ns}t_{j}3^{-j}$,
$t_{j}\in\cbr{0,1,2}$, $j=1,2,\ldots$, be the trinomial expansion
of the number $t\in[0,1]$. Let $k:\cbr{0,1,2}\ra\cbr{0,1,2}$ be
given by $k(n)=2-n$. Now, the Peano curve ${\cal P}(t)=(x(t),y(t))$
is defined by the formulas
\begin{equation}
x(t)=0_{3} k^{0}\rbr{t_{1}}k^{t_{2}}\rbr{t_{3}}k^{t_{2}+t_{4}}\rbr{t_{5}}k^{t_{2}+t_{4}+t_{6}}\rbr{t_{7}}\ldots,\label{eq:px}
\end{equation}
\begin{equation}
y(t)=0_{3} k^{t_{1}}\rbr{t_{2}}k^{t_{1}+t_{3}}\rbr{t_{4}}k^{t_{1}+t_{3}+t_{5}}\rbr{t_{6}}\ldots;\label{eq:py}
\end{equation}
here $k^{v}$, $v=0,1,2,\ldots$ denotes $v$th iteration of the function
$k$. For the proof that formulas (\ref{eq:px}) and (\ref{eq:py})
define $x$ and $y$ in a unique way (regardless whether we take finite
or infinite trinomial expansion of $t$ (with infinite many trailing
$2$s) \textendash{} in case it has finite trinomial expansion), as
well as that ${\cal P}$ is continuous and the image of the segment $[0,1]$
under ${\cal P}$ is the whole unit square $[0,1]^{2}$ see for example
\cite[Chapt. 3]{Sagan:1994}. For illustration, we present
two approximating polygons of the Peano curve obtained in the first and the second step of an iteration procedure, whose limit is the Peano curve. 
\begin{center}
\includegraphics[scale=0.65]{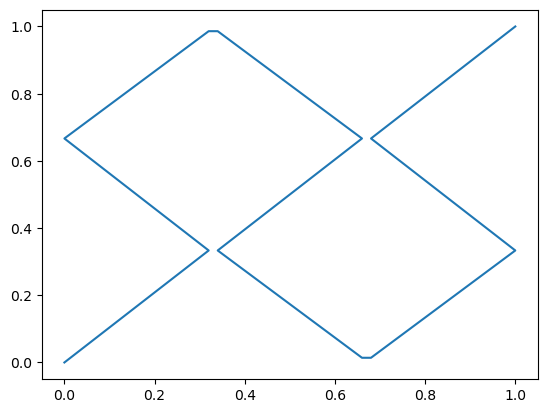}
\captionof{figure}{Approximating polygon of the Peano curve (first step of iteration procedure). Some lines were altered for better illustration of the parametrization of the curve.}
\includegraphics[scale=0.65]{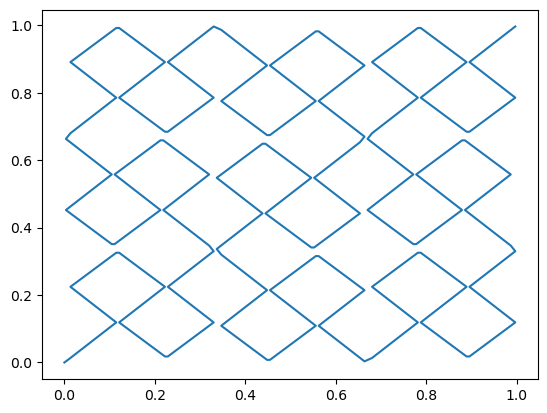}
\captionof{figure}{Approximating polygon of the Peano curve (second step of iteration procedure). Some lines were altered for better illustration of the iteration procedure.}
\end{center}
Let $P$ be the approximating polygon obtained in the first step. In the second step each segment of $P$ is replaced by the whole $P$ (after proper rescaling, shift and/or reflection). Similarly, in the $n$th step each segment of approximating polygon obtained in the previous step is replaced by the whole $P$ (after proper rescaling, shift and/or reflection).

The $x$ component of the first two approximating polygons of
the Peano curve is presented in the two next figures.
\begin{center}
\includegraphics[scale=0.65]{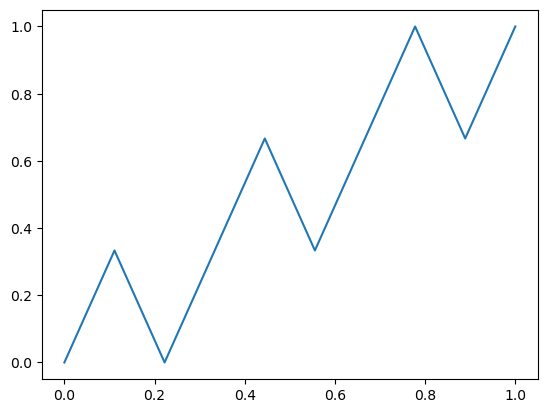}
\captionof{figure}{The $x$ component of the first approximating polygon.}
\end{center}
\begin{center}
\includegraphics[scale=0.65]{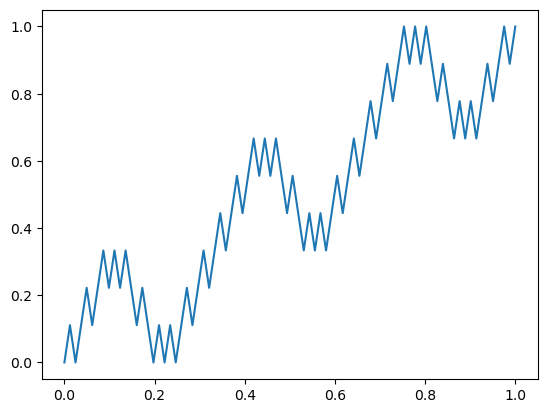}
\captionof{figure}{The $x$ component of the second approximating polygon.}
\end{center}

Using (\ref{eq:px}) one easily checks that $x:[0,1]\ra[0,1]$
is self-similar in this regard that
\begin{equation}
x(t)=\frac{1}{3}x(9t),\quad t\in\sbr{0,\frac{1}{9}},\label{eq:s1}
\end{equation}
\begin{equation}
x(t)=\frac{1}{3}x(2-9t) = x\rbr{\frac{2}{9}-t},\quad t\in\sbr{\frac{1}{9},\frac{2}{9}},\label{eq:s2}
\end{equation}
\begin{equation}
x(t)=\frac{1}{3}x(9t-2) = x\rbr{t-\frac{2}{9}},\quad t\in\sbr{\frac{2}{9},\frac{3}{9}},\label{eq:s3}
\end{equation}
\begin{equation}
x(t)=\frac{1}{3}+\frac{1}{3}x(9t-3) = \frac{1}{3}+x\rbr{t-\frac{3}{9}},\quad t\in\sbr{\frac{3}{9},\frac{4}{9}},\label{eq:s4}
\end{equation}
\begin{equation}
x(t)=\frac{1}{3}+\frac{1}{3}x(5-9t)=  \frac{1}{3}+x\rbr{\frac{5}{9}-t}, \quad t\in\sbr{\frac{4}{9},\frac{5}{9}},\label{eq:s5}
\end{equation}
\begin{equation}
x(t)=\frac{1}{3}+\frac{1}{3}x(9t-5) = \frac{1}{3}+x\rbr{t-\frac{5}{9}}, \quad t\in\sbr{\frac{5}{9},\frac{6}{9}},\label{eq:s6}
\end{equation}
\begin{equation}
x(t)=\frac{2}{3}+\frac{1}{3}x(9t-6) = \frac{2}{3}+x\rbr{t-\frac{6}{9}},\quad t\in\sbr{\frac{6}{9},\frac{7}{9}},\label{eq:s7}
\end{equation}
\begin{equation}
x(t)=\frac{2}{3}+\frac{1}{3}x(8-9t) = \frac{2}{3}+x\rbr{\frac{8}{9}-t},\quad t\in\sbr{\frac{7}{9},\frac{8}{9}},\label{eq:s8}
\end{equation}
\begin{equation}
x(t)=\frac{2}{3}+\frac{1}{3}x(9t-8) = \frac{2}{3}+x\rbr{t-\frac{8}{9}}, \quad t\in\sbr{\frac{8}{9},\frac{9}{9}}.\label{eq:s9}
\end{equation}
We also have the following symmetry 
\begin{equation}
x(1-t)=1-x(t), \quad t\in\sbr{0,1}.\label{eq:s10}
\end{equation}

A very similar self-similarity holds for the standard Brownian motion
stochastic process $\rbr{B(t),t\in[0,+\ns)}$, namely
\[
\rbr{B(t),\quad t\in\sbr{0,\frac{1}{9}}}=^{d}\rbr{\frac{1}{3}B(9t),\quad t\in\sbr{0,\frac{1}{9}}},
\]
where '$=^{d}$' denotes the equality of finite dimensional distributions
of processes $\rbr{B(t),\quad t\in\sbr{0,\frac{1}{9}}}$ and $\rbr{3^{-1}B(9t),\quad t\in\sbr{0,\frac{1}{9}}}$.
Arguably, a similar but more general property of the standard Brownian
motion (sBm in short), namely
\[
\rbr{B(t)-B(s),\quad0\le s\le t}=^{d}\rbr{a^{-1}\rbr{B(a^{2}t)-B(a^{2}s)},\quad0\le s\le t},\text{ for any }a>0,
\]
is the reason why the sBm has linear quadratic variation, $[B]_t = t, t\ge 0$. Quadratic variation $[B]_t$, $ t\ge 0$, of the sBm is usually defined as the limit of the sums
\[
\sum_{i=1}^{k_n} \rbr{B\rbr{t_i^n} - B\rbr{t_{i-1}^n}}^2,
\]
called \emph{quadratic variations along (or based on) the sequences $\rbr{t_i^n}_{i = 0}^{k_n} $}, $n=1,2,\ldots$,  where $$0=t_0^n < t_1^n < t_2^n < \ldots < t_{k_n}^n = t$$ are deterministic times such that $$\max_{i=1,2,\ldots,k_n} \rbr{t_i^n - t_{i-1}^n}$$ tends to $0$ as $n \ra +\ns$. For sBm, this limit exists (in probability) and is equal $t$.

A better, almost sure convergence, is obtained when $t_i^n$s are replaced by \emph{random} times, which are times of the hitting by $B$ the grid $c_n \Z+r_n$, where $c_n$ is a sequence of positive reals tending to $0$, $r_n \in \R$, see \cite{Chacon:1981}. Partitions obtained from the hitting times by a path the grids $c_n \Z + r_n$ are called \emph{Lebesgue partitions}; for a precise definition in a case where a path is continuous see Sect. \ref{LebPar}.
A similar phenomenon holds for quadratic variation of continuous semimartingales, see \cite{Lemieux:1983}, c\`adl\`ag semimartingales \cite{LochOblPS:2021}, and $1/H$-variation of fractional Brownian motions (fBms in short) \cite{Toyomu:2023}.

One may expect a similar phenomenon for the $x$-component of the Peano curve. It is not difficult to observe that quadratic variations of $x$ based on Lebesgue partitions for the hitting times by $x$ the grids $3^{-n} \Z + r_n$, where $r_n$ is not multiple of $3^{-n}$, $n=1,2,\ldots $, tend to $t/4$. However, computer simulations (a Python language code available upon request) revealed quite different picture for other grids of the form $3^{-n} p \Z + r_n$, where $p$ is a rational number. In the table below we present some results of these simulations.
\begin{center}
\begin{tabular}{ |c|c|c| } 
\hline
 $c_n =$ & $r_n =$ & simulation of the quadratic variation \\ 
 \hline
 \hline
 $1/3^{n}$ & $0$ & $ t$ \\ 
  $2/3^{n}$ & $0$ & $ 0.5 \ t$ \\ 
  $1/3^{n}$ & $1/\rbr{2\cdot 3^n}$ & $ 0.25 \ t$ \\ 
  $2/ 3^n$ & $1/\rbr{2 \cdot 3^n}$ & $0.3333 \ t$ \\   
  $2/ 3^n$ & $1/\rbr{4 \cdot 3^n}$ & $0.3333 \ t$ \\ 
 $1/\rbr{2 \cdot 3^n}$ & $1/\rbr{4 \cdot 3^n}$ & $0.3125 \ t$ \\ 
$8/ 3^n$ & $1/\rbr{5 \cdot 3^n}$ & $0.2963 \ t$ \\ 
\hline
\end{tabular}
\end{center}
Thus we observe that the  situation is quite different for the Peano curve -- the quadratic
variations of the $x$-component of the Peano curve along different sequences
of Lebesgue partitions may be quite different (or may even not exist). The pattern emerging from the simulation results is not clear. However, we managed to obtain, using rather delicate and detailed reasoning, the existence of the   quadratic
variations of the $x$-component of the Peano curve along some sequences of Lebesgue partitions.
More precisely, we proved that the $x$-component of the Peano curve has quadratic
variation defined as the limit (as $n \ra +\ns$) of quadratic variations along the Lebesgue
partitions for grids of the form $3^{-n}p\Z+3^{-n}r$, $n=1,2,\ldots$, where $p$ is a positive rational number and $r$ is for example an irrational number, and this quadratic variation is given by
\[
[x]_t = 3^{\left\lfloor -\log_{3}p\right\rfloor }p\cbr{1-\frac{3}{4}\rbr{3^{\left\lfloor -\log_{3}p\right\rfloor }p}}t, \quad t \in [0,1].
\]
This stands in contrast with the property of the trajectories of above mentioned
stochastic processes (independence of the quadratic variation along Lebesgue partitions from the grids) and reflects \emph{regularisation by noise phenomenon}.

More detailed information about trajectories of a stochastic process
than the quadratic (or $1/H$) variation is provided by \emph{local
time}. Local time is usually defined as the density of the occupation
measure of the trajectories of the process with respect to some given
measure $\nu$. For fBms, the measure $\nu$ is the Lebesgue measure,
for semimartingales - the measure $\nu$ is the Lebesgue-Stieltjes
measure associated with the continuous part of the quadratic variation of the process. Such
defined local time for fBms and semimartingales may be defined in
terms of grid crossings, but now we consider only points of the grid
around the given level.

In what follows, we prove that the Peano curve has local time defined
as the density of its occupation measure with respect to the Lebesgue
measure, and it can be expressed in terms of numbers of grid crossings around
the given level, but the normalization of these numbers is not as smooth as in the case of stochastic processes. This is another manifestation of regularisation by noise phenomenon.

{The paper is organized as follows. In Section 1.1 we state essential definitions used later in the paper. Section 2 explores the quadratic variation of the Peano curve along the sequence of Lebesgue partitions. In Section 2.1, we introduce a one-step recursion for the quadratic variation of the $x$-component of the Peano curve along these partitions, which is then extended to a $k$-step recursion in Section 2.2. Section 2.3 examines the limit of quadratic variations of $x$ along the Lebesgue partitions and presents the first main result of the paper. In Section 2.4, we discuss truncated variation and the number of crossings of the $x$-component. In section 3 we prove the second main result of this paper -- existence of the local times of the $x$-component of the Peano curve, and expresses this local time as the weak limit of normalized interval crossing counts.}

Notation: $\Z$ denotes the set of integers, $\N$ denotes the set of non-negative integers, 
$\N_{+}$ denotes the set of positive integers and $\Q$ denotes the set of rational numbers.

\subsection{Quadratic variation of a continuous curve based on the Lebesgue partition -- definitions} \label{LebPar}

Let $c$ be a positive real and $r$ be a real number. By the grid
$c\Z+r$ we will mean the countable set
\[
c\Z+r=\cbr{u\in\R:\exists q\in\Z\quad u=c\cdot q+r}.
\]

Now let us define the Lebesgue partition $\pi(c,r)$ for a continuous
function $x:[0,T]\ra\R$ ($T$ is a positive real) and the grid $c\Z+r$. $\pi(c,r)$ is the
sequence of consecutive points $t_{1}(c,r)\le t_{2}(c,r)\le\ldots$ from
the set $[0,T]\cup\cbr{+\ns}$,
\[
\pi(c,r)=\rbr{t_{1}(c,r),t_{2}(c,r),\ldots},
\]
defined recursively in the following way.
\[
t_{1}(c,r)=\inf\cbr{t\in[0,T]:x(t)\in\rbr{c\Z+r}}
\]
with the convention that $\inf\emptyset=+\ns$. If $t_{m}(c,r)=+\ns$
for some $m=1,2,\ldots$, then we set $t_{m+1}(c,r)=+\ns$. If $t_{m}(c,r)<+\ns$
for some $m=1,2,\ldots$, (this also means that $t_{m}(c,r)\le T$)
then we set
\[
t_{m+1}(c,r)=\inf\cbr{t\in\sbr{t_{m}(c,r),T}:x(t)\in\rbr{\rbr{c\Z+r}\setminus\cbr{x\rbr{t_{m}(c,r)}}}}.
\]

For $t\in[0,T]$ we also define
\[
k(c,r,t):=\min\cbr{m\in\N_{+}:t_{m}(c,r)\ge t},\quad l(c,r,t):=\max\cbr{m\in\N_{+}:t_{m}(c,r)\le t}
\]
and the quadratic variation $\sbr x_{s,t}^{c,r}$ of $x:[0,T]\ra\R$
on the time interval $[s,t]\subseteq\sbr{0,T}$ along the Lebesgue
partition $\pi\rbr{c,r}$ for $x$ is defined as
\[
\sbr x_{s,t}^{c,r}:=\sum_{m=k(c,r,s)}^{l(c,r,t)-1}\rbr{x\rbr{t_{m+1}(c,r)}-x\rbr{t_{m}(c,r)}}^{2},
\]
where the value of the sum $\sum_{m=k(c,r,s)}^{l(c,r,t)-1}$ is $0$
when $l(c,r,s)-1<k(c,r,s)$. Let us notice that since $\left|x\rbr{t_{m+1}(c,r)}-x\rbr{t_{m}(c,r)}\right|=c$
for $k(c,r,t)\le m\le l(c,r,t)-1$ then we have
\begin{equation}
\sbr x_{s,t}^{c,r}=\rbr{l(c,r,t)-k(c,r,t)}c^{2}.\label{eq:qvar_leb}
\end{equation}
We also have the following relation: for $s<t<u$
\begin{equation}
\sbr x_{s,u}^{c,r}=\sbr x_{s,t}^{c,r}+{\bf 1}_{l(c,r,t)<k(c,r,t)}c^{2}+\sbr x_{t,u}^{c,r},\label{eq:intervals_quasi_additivity}
\end{equation}
where
\[
{\bf 1}_{l(c,r,t)<k(c,r,t)}=\begin{cases}
0 & \text{ if }l(c,r,t)\ge k(c,r,t),\\
1 & \text{ if }l(c,r,t)<k(c,r,t)
\end{cases}=\begin{cases}
0 & \text{ if }l(c,r,t)=k(c,r,t)=t,\\
1 & \text{ if }l(c,r,t)<t<k(c,r,t).
\end{cases}
\]

For $t\in\sbr{0,T}$, by $\sbr x_{t}^{c,r}$we will mean $\sbr x_{0,t}^{c,r}$,
so we formally define \textbf{$\sbr x_{t}^{c,r}$} as
\[
\sbr x_{t}^{c,r}:=\sbr x_{0,t}^{c,r}.
\]

\section{Quadratic variation of the Peano curve along the sequence of the
Lebesgue partitions for the grids $\protect\rbr{p/3^{n}}\protect\Z+\protect\rbr{r/3^{n}}$,
$p\in\protect\Q$, $n\in\protect\N$}

Let $x:[0,1]\ra\R$ be the $x$ component of the Peano curve.
In this section we will deal with the quadratic variations of $x$ along the Lebesgue
partitions $\pi\rbr{c_n,r_n}$ for $x$ and the grids $c_n\Z+r_n$ where
 \[
c_{n}=\frac{p}{3^{n}} \le \frac{1}{3},\quad  r =\frac{r}{3^{n}},
\]
$p=p'/q'$, $p'$ and $q'$ are positive integers not divisible
by $3$, $r\in \R$, and $n\in\Z$. 
\subsection{One-step recursion}

Let us fix $c\in(0, 1/3]$ and $r\in \R$. From self-similarity relation (\ref{eq:s1}), we see that the range
of the function $x$ on the interval $[0,1/9]$ is $[0,1/3]$ and
\[
t_{m}\rbr{c,r}=\frac{1}{9}t_{m}\rbr{3c,3r},\quad k\rbr{c,r,1/9}=k\rbr{3c,3r,1},\quad l\rbr{c,r,1/9}=l\rbr{3c,3r,1}.
\]
This and (\ref{eq:qvar_leb}) yield
\begin{equation}
\sbr x_{1/9}^{c,r}=\frac{1}{9}\sbr x_{1}^{3c,3r}\label{eq:019}
\end{equation}
(recall that $\sbr x_{1/9}^{c,r}=\sbr x_{0,1/9}^{c,r}$, $\sbr x_{1}^{3c,3r}=\sbr x_{0,1}^{3c,3r}$).
Similarly, by (\ref{eq:s2}) and using the fact that the points from
the grid $c\Z+r$ hit by $x$ consecutively on the time interval $\sbr{1/9,2/9}$
are the same as the points from this grid hit by $x$ on the time
interval $\sbr{0,1/9}$ but in the reverse order, we get
\begin{equation}
\sbr x_{2/9}^{c,r}-\sbr x_{1/9}^{c,r}=\sbr x_{1/9}^{c,r}=\frac{1}{9}\sbr x_{1}^{3c,3r}.\label{eq:192939}
\end{equation}
Using (\ref{eq:s3}) and a similar observation for the intervals $\sbr{1/9,2/9}$
and $\sbr{2/9,3/9}$ we get
\begin{equation}
\sbr x_{3/9}^{c,r}-\sbr x_{2/9}^{c,r}=\sbr x_{1/9}^{c,r}=\frac{1}{9}\sbr x_{1}^{3c,3r}.\label{eq:192939-1}
\end{equation}

Next, we will use relations (\ref{eq:s4})-(\ref{eq:s6}), or equivalent
relations
\begin{equation}
3x(t)-1=x(9t-3),\quad t\in\sbr{\frac{3}{9},\frac{4}{9}},\label{eq:s4-1}
\end{equation}
\begin{equation}
3x(t)-1=x(5-9t),\quad t\in\sbr{\frac{4}{9},\frac{5}{9}},\label{eq:s5-1}
\end{equation}
\begin{equation}
3x(t)-1=3x(9t-5),\quad t\in\sbr{\frac{5}{9},\frac{6}{9}}.\label{eq:s6-1}
\end{equation}
From the fact that the function $x$ attains on the interval $\sbr{3/9,4/9}$
values from the interval $[1/3,2/3]$, while on the interval $\sbr{0,1/3}$
it attains values from the interval $[0,1/3]$ it follows that either
the first value from the grid $c\Z+r$ attained by $x$ on the time
interval $[3/9,4/9]$ is different from the last value from the grid
$c\Z+r$ attained by $x$ on the time interval $[0,1/3]$, or this
value is equal $1/3$ and is attained by $x$ at $t=1/3$. This together
with relations (\ref{eq:s4-1}) and (\ref{eq:intervals_quasi_additivity})
yield

\begin{equation}
\sbr x_{3/9,4/9}^{c,r}=\sbr x_{4/9}^{c,r}-\sbr x_{1/3}^{c,r}-{\bf 1}_{l(c,r,1/3)<k(c,r,1/3)}c^{2}=\frac{1}{9}\sbr x_{1}^{3c,3r-1}.\label{eq:3949}
\end{equation}
Next, (\ref{eq:s5-1}), (\ref{eq:s6-1}), the fact that the points
from the grid $c\Z+r$ hit by $x$ consecutively on the time interval
$\sbr{4/9,5/9}$ are the same as the points from this grid hit by
$x$ on the time interval $\sbr{3/9,4/9}$ but in the reverse order,
and a similar observation for the intervals $\sbr{4/9,5/9}$ and $\sbr{5/9,6/9}$
give
\begin{equation}
\sbr x_{5/9}^{c,r}-\sbr x_{4/9}^{c,r}=\sbr x_{3/9,4/9}^{c,r}=\frac{1}{9}\sbr x_{1}^{3c,3r-1},\label{eq:4959}
\end{equation}
\begin{equation}
\sbr x_{6/9}^{c,r}-\sbr x_{5/9}^{c,r}=\sbr x_{3/9,4/9}^{c,r}=\frac{1}{9}\sbr x_{1}^{3c,3r-1}.\label{eq:5969}
\end{equation}

Similarly, for the intervals $\sbr{6/9,7/9}$, $\sbr{7/9,8/9}$ and
$\sbr{8/9,9/9}$ one has, respectively,
\begin{equation}
\sbr x_{6/9,7/9}^{c,r}=\sbr x_{7/9}^{c,r}-\sbr x_{2/3}^{c,r}-{\bf 1}_{l(c,r,2/3)<k(c,r,2/3)}c^{2}=\frac{1}{9}\sbr x_{1}^{3c,3r-2},\label{eq:6979}
\end{equation}
\begin{equation}
\sbr x_{8/9}^{c,r}-\sbr x_{7/9}^{c,r}=\sbr x_{6/9,7/9}^{c,r}=\frac{1}{9}\sbr x_{1}^{3c,3r-2}\label{eq:7989}
\end{equation}
and
\begin{equation}
\sbr x_{9/9}^{c,r}-\sbr x_{8/9}^{c,r}=\sbr x_{6/9,7/9}^{c,r}=\frac{1}{9}\sbr x_{1}^{3c,3r-2}.\label{eq:8999}
\end{equation}

Summing equalities (\ref{eq:019})-(\ref{eq:192939-1}), (\ref{eq:3949})-(\ref{eq:5969})
and (\ref{eq:6979})-(\ref{eq:8999}) we get
\[
\sbr x_{1}^{c,r}-{\bf 1}_{l(c,r,1/3)<k(c,r,1/3)}c^{2}-{\bf 1}_{l(c,r,2/3)<k(c,r,2/3)}c^{2}=\frac{1}{3}\sbr x_{1}^{3c,3r}+\frac{1}{3}\sbr x_{1}^{3c,3r-1}+\frac{1}{3}\sbr x_{1}^{3c,3r-2}
\]
or, equivalently
\begin{equation}
\sbr x_{1}^{c,r}=\frac{1}{3}\sbr x_{1}^{3c,3r}+\frac{1}{3}\sbr x_{1}^{3c,3r-1}+{\bf 1}_{l(c,r,1/3)<k(c,r,1/3)}c^{2}+\frac{1}{3}\sbr x_{1}^{3c,3r-2}+{\bf 1}_{l(c,r,2/3)<k(c,r,2/3)}c^{2}.\label{eq:recccursion}
\end{equation}

\subsection{$k$-step recursion}

Now we will deal with the quadratic variations of the function $x$
along the Lebesgue partitions $\pi\rbr{c_n,r_n}$ for $x$ and the grids
$c_{n}\Z+r_{n}$ where
\[
c_{n}\Z+r_{n}=\frac{p}{3^{n}}\Z+\frac{r}{3^{n}}=\frac{p'}{q'\cdot3^{n}}\Z+\frac{r}{3^{n}},
\]
$p=p'/q'$, $p'$ and $q'$ are positive integers not divisible
by $3$, $r\in \R$, and $n\in\Z$. This means that
\[
c_{n}=\frac{p}{3^{n}}=\frac{p'}{3^{n}q'},\quad r_{n}=\frac{r}{3^{n}},\quad n\in\Z.
\]
Let us define
\begin{equation}
\mathbb{M}:=\cbr{\tilde{n}\in\Z:\frac{p}{3^{\tilde{n}}}\le\frac{1}{3}}.\label{eq:Mdef}
\end{equation}
Introducing the new variable $\theta=r/p=rq'/p'\in[0,1)$ such that
$r=\theta c$ and using the equality $3c_{n}=c_{n-1}$ for $n\in \M$,
by (\ref{eq:recccursion}) we get
\begin{align}
\sbr x_{1}^{c_{n},r_{n}} & =\frac{1}{3}\sbr x_{1}^{c_{n-1},\theta c_{n-1}}+\frac{1}{3}\sbr x_{1}^{c_{n-1},\theta c_{n-1}-1}+{\bf 1}_{l(c_{n},r_{n},1/3)<k(c_{n},r_{n},1/3)}c^{2}\label{eq:rekursja0}\\
 & \quad+\frac{1}{3}\sbr x_{1}^{c_{n-1},\theta c_{n-1}-2}+{\bf 1}_{l(c_{n},r_{n},2/3)<k(c_{n},r_{n},2/3)}c^{2}.
\end{align}
Let us notice that the inequality
\[
l(c_{n},r_{n},1/3)<k(c_{n},r_{n},1/3)
\]
holds iff $l(c_{n},r_{n},1/3)<1/3<k(c_{n},r_{n},1/3)$ which holds
iff the grid $c_{n}\Z+r_{n}$ does not contain the value $x(1/3)=1/3$
which is equivalent with the statements
\[
\forall q\in\Z\quad c_{n}\cdot q+r_{n}=c_{n}\rbr{q+\theta}\neq\frac{1}{3},
\]
\[
\forall q\in\Z\quad q+\theta\neq\frac{1}{c_{n-1}},
\]
\[
\cbr{\theta-\frac{1}{c_{n-1}}}\neq0,
\]
where $\cbr a$ denotes the fractional part of the number $a$. Similarly,
the inequality
\[
l(c_{n},r_{n},2/3)<k(c_{n},r_{n},2/3)
\]
holds iff
\[
\left\{ \theta-\frac{2}{c_{n-1}}\right\} \neq0.
\]
Now, (\ref{eq:rekursja0}) takes for $n\in\M$ the form
\begin{align}
\sbr x_{1}^{c_{n},\theta c_{n}} & =\frac{1}{3}\sbr x_{1}^{c_{n-1},\theta c_{n-1}}\nonumber \\
 & \quad+\frac{1}{3}\sbr x_{1}^{c_{n-1},\left\{ \theta-1/c_{n-1}\right\} c_{n-1}}+\mathbf{1}_{(0,1)}\rbr{\left\{ \theta-\frac{1}{c_{n-1}}\right\} }c^{2}\nonumber \\
 & \quad+\frac{1}{3}\sbr x_{1}^{c_{n-1},\left\{ \theta-2/c_{n-1}\right\} c_{n-1}}+\mathbf{1}_{(0,1)}\rbr{\left\{ \theta-\frac{2}{c_{n-1}}\right\} }c^{2},\label{eq:recursion_impr}
\end{align}
where $\mathbf{1}_{(0,1)}$ is the indicator function of the interval
$(0,1)$:
\[
\mathbf{1}_{(0,1)}(a)=\begin{cases}
0 & \text{ if }a\notin(0,1),\\
1 & \text{ if }a\in(0,1).
\end{cases}
\]
By (\ref{eq:recursion_impr}) we have the following $1$-step 
recursion

\begin{align}
\sbr x_{1}^{c_{n},\theta c_{n}} & =\frac{1}{3}\sbr x_{1}^{c_{n-1},\theta c_{n-1}}\nonumber \\
 & \quad+\frac{1}{3}\sbr x_{1}^{c_{n-1},\left\{ \theta-1/c_{n-1}\right\} c_{n-1}}+\mathbf{1}_{(0,1)}\rbr{\left\{ \theta-\frac{1}{c_{n-1}}\right\} }c^{2}\nonumber \\
 & \quad+\frac{1}{3}\sbr x_{1}^{c_{n-1},\left\{ \theta-2/c_{n-1}\right\} c_{n-1}}+\mathbf{1}_{(0,1)}\rbr{\left\{ \theta-\frac{2}{c_{n-1}}\right\} }c^{2}\nonumber \\
 & =\frac{1}{3}\sbr x_{1}^{\frac{3p'}{3^{n}q'},\theta\frac{3p'}{3^{n}q'}}\nonumber \\
 & \quad+\frac{1}{3}\sbr x_{1}^{\frac{3p'}{3^{n}q'},\left\{ \theta-\frac{3^{n}q'}{3p'}\right\} \frac{3p'}{3^{n}q'}}+\mathbf{1}_{(0,1)}\rbr{\left\{ \theta-\frac{3^{n}q'}{3p'}\right\} }\rbr{\frac{p'}{3^{n}q'}}^{2}\nonumber \\
 & \quad+\frac{1}{3}\sbr x_{1}^{\frac{3p'}{3^{n}q'},\left\{ \theta-2\frac{3^{n}q'}{3p'}\right\} \frac{3p'}{3^{n}q'}}+\mathbf{1}_{(0,1)}\rbr{\left\{ \theta-2\frac{3^{n}q'}{3p'}\right\} }\rbr{\frac{p'}{3^{n}q'}}^{2}.\label{eq:one_step-1}
\end{align}

Now, let us assume that
\[
n-1\in\M
\]
which means that
\[
\frac{p'}{3^{n-1}q'}\le\frac{1}{3},\text{ so }\frac{3^{2}p'}{3^{n}q'}\le1
\]
and we can do the second step:
\begin{align*}
\sbr x_{1}^{\frac{3p'}{3^{n}q'},\theta\frac{3p'}{3^{n}q'}}=& \frac{1}{3} \sbr x_{1}^{\frac{3^{2}p'}{3^{n}q'},\theta\frac{3^{2}p'}{3^{n}q'}}\\
 &+\frac{1}{3}\sbr x_{1}^{\frac{3^{2}p'}{3^{n}q'},\left\{ \theta-\frac{3^{n}q'}{3^{2}p'}\right\} \frac{3^{2}p'}{3^{n}q'}}+\mathbf{1}_{(0,1)}\rbr{\left\{ \theta-\frac{3^{n}q'}{3^{2}p'}\right\} }\rbr{\frac{3p'}{3^{n}q'}}^{2}\\
 & +\frac{1}{3}\sbr x_{1}^{\frac{3^{2}p'}{3^{n}q'},\left\{ \theta-2\frac{3^{n}q'}{3^{2}p'}\right\} \frac{3^{2}p'}{3^{n}q'}}+\mathbf{1}_{(0,1)}\rbr{\left\{ \theta-2\frac{3^{n}q'}{3^{2}p'}\right\} }\rbr{\frac{3p'}{3^{n}q'}}^{2}.
\end{align*}
Similarly we can expand two other terms appearing in  \eqref{eq:one_step-1}:
\begin{align*}
\sbr x_{1}^{\frac{3p'}{3^{n}q'},\left\{ \theta-\frac{3^{n}q'}{3p'}\right\} \frac{3p'}{3^{n}q'}}= & \frac{1}{3} \sbr x_{1}^{\frac{3^{2}p'}{3^{n}q'},\left\{ \theta-\frac{3^{n}q'}{3p'}\right\} \frac{3^{2}p'}{3^{n}q'}}\\
 & +\frac{1}{3}\sbr x_{1}^{\frac{3^{2}p'}{3^{n}q'},\left\{ \theta-\frac{3^{n}q'}{3p'}-\frac{3^{n}q'}{3^{2}p'}\right\} \frac{3^{2}p'}{3^{n}q'}}+\mathbf{1}_{(0,1)}\rbr{\left\{ \theta-\frac{3^{n}q'}{3p'}-\frac{3^{n}q'}{3^{2}p'}\right\} }\rbr{\frac{3p'}{3^{n}q'}}^{2}\\
 & +\frac{1}{3}\sbr x_{1}^{\frac{3^{2}p'}{3^{n}q'},\left\{ \theta-\frac{3^{n}q'}{3p'}-2\frac{3^{n}q'}{3^{2}p'}\right\} \frac{3^{2}p'}{3^{n}q'}}+\mathbf{1}_{(0,1)}\rbr{\left\{ \theta-\frac{3^{n}q'}{3p'}-2\frac{3^{n}q'}{3^{2}p'}\right\} }\rbr{\frac{3p'}{3^{n}q'}}^{2}
\end{align*}
and
\begin{align*}
\sbr x_{1}^{\frac{3p'}{3^{n}q'},\left\{ \theta-2\frac{3^{n}q'}{3p'}\right\} \frac{3p'}{3^{n}q'}}= & \frac{1}{3}\sbr x_{1}^{\frac{3^{2}p'}{3^{n}q'},\left\{ \theta-2\frac{3^{n}q'}{3p'}\right\} \frac{3^{2}p'}{3^{n}q'}}\\
 & +\frac{1}{3}\sbr x_{1}^{\frac{3^{2}p'}{3^{n}q'},\left\{ \theta-2\frac{3^{n}q'}{3p'}-\frac{3^{n}q'}{3^{2}p'}\right\} \frac{3^{2}p'}{3^{n}q'}}+\mathbf{1}_{(0,1)}\rbr{\left\{ \theta-2\frac{3^{n}q'}{3p'}-\frac{3^{n}q'}{3^{2}p'}\right\} }\rbr{\frac{3p'}{3^{n}q'}}^{2}\\
 & +\frac{1}{3}\sbr x_{1}^{\frac{3^{2}p'}{3^{n}q'},\left\{ \theta-2\frac{3^{n}q'}{3p'}-2\frac{3^{n}q'}{3^{2}p'}\right\} \frac{3^{2}p'}{3^{n}q'}}+\mathbf{1}_{(0,1)}\rbr{\left\{ \theta-2\frac{3^{n}q'}{3p'}-2\frac{3^{n}q'}{3^{2}p'}\right\} }\rbr{\frac{3p'}{3^{n}q'}}^{2}.
\end{align*}
From the first and the second steps we get
\begin{align*}
\sbr x_{1}^{c_{n},\theta c_{n}} & =\frac{1}{9}\sum_{i_{1},i_{2}\in\cbr{0,1,2}}\sbr x_{1}^{\frac{3^{2}p'}{3^{n}q'},\left\{ \theta-i_{1}\frac{3^{n}q'}{3p'}-i_{2}\frac{3^{n}q'}{3^{2}p'}\right\} \frac{3^{2}p'}{3^{n}q'}}\\
 & \quad+\mathbf{1}_{(0,1)}\rbr{\left\{ \theta-\frac{3^{n}q'}{3p'}\right\} }\rbr{\frac{p'}{3^{n}q'}}^{2}+\mathbf{1}_{(0,1)}\rbr{\left\{ \theta-2\frac{3^{n}q'}{3p'}\right\} }\rbr{\frac{p'}{3^{n}q'}}^{2}\\
 & \quad+\frac{1}{3}\mathbf{1}_{(0,1)}\rbr{\left\{ \theta-\frac{3^{n}q'}{3^{2}p'}\right\} }\rbr{\frac{3p'}{3^{n}q'}}^{2}+\frac{1}{3}\mathbf{1}_{(0,1)}\rbr{\left\{ \theta-2\frac{3^{n}q'}{3^{2}p'}\right\} }\rbr{\frac{3p'}{3^{n}q'}}^{2}\\
 & \quad+\frac{1}{3}\mathbf{1}_{(0,1)}\rbr{\left\{ \theta-\frac{3^{n}q'}{3p'}-\frac{3^{n}q'}{3^{2}p'}\right\} }\rbr{\frac{3p'}{3^{n}q'}}^{2}+\frac{1}{3}\mathbf{1}_{(0,1)}\rbr{\left\{ \theta-\frac{3^{n}q'}{3p'}-2\frac{3^{n}q'}{3^{2}p'}\right\} }\rbr{\frac{3p'}{3^{n}q'}}^{2}\\
 & \quad+\frac{1}{3}\mathbf{1}_{(0,1)}\rbr{\left\{ \theta-2\frac{3^{n}q'}{3p'}-\frac{3^{n}q'}{3^{2}p'}\right\} }\rbr{\frac{3p'}{3^{n}q'}}^{2}+\frac{1}{3}\mathbf{1}_{(0,1)}\rbr{\left\{ \theta-2\frac{3^{n}q'}{3p'}-2\frac{3^{n}q'}{3^{2}p'}\right\} }\rbr{\frac{3p'}{3^{n}q'}}^{2}\\
 & =\frac{1}{3^{2}}\sum_{i_{1},i_{2}\in\cbr{0,1,2}}\sbr x_{1+}^{\frac{3^{2}p'}{3^{n}q'},\left\{ \theta-\sum_{v=1}^{2}i_{v}\frac{3^{n}q'}{3^{v}p'}\right\} \frac{3^{2}p'}{3^{n}q'}}\\
 & \quad+\sum_{m=1}^{2}\frac{1}{3^{m-1}}\sum_{j_{1},\ldots,j_{m-1}\in\cbr{0,1,2},j_{m}\in\cbr{1,2}}\mathbf{1}_{(0,1)}\rbr{\left\{ \theta-\sum_{l=1}^{m}j_{l}\frac{3^{n}q'}{3^{l}p'}\right\} }\rbr{\frac{3^{m-1}p'}{3^{n}q'}}^{2}.
\end{align*}

Reasoning similarly further, for a given integer $n$ and any positive
integer $k$ such that $n-(k-1)\in\M$, we obtain
\begin{align}
\sbr x_{1}^{c_{n},\theta c_{n}} & =\frac{1}{3^{k}}\sum_{i_{1},i_{2},\ldots,i_{k}\in\cbr{0,1,2}}\sbr x_{1}^{\frac{3^{k}p'}{3^{n}q'},\left\{ \theta-\sum_{v=1}^{k}i_{v}\frac{3^{n}q'}{3^{v}p'}\right\} \frac{3^{k}p'}{3^{n}q'}}\label{eq:k-step_rec}\\
 & \quad+\sum_{m=1}^{k}\frac{1}{3^{m-1}}\sum_{j_{1},...,j_{m-1}\in\cbr{0,1,2},j_{m}\in\cbr{1,2}}\mathbf{1}_{(0,1)}\rbr{\left\{ \theta-\sum_{l=1}^{m}j_{l}\frac{3^{n}q'}{3^{l}p'}\right\} }\rbr{\frac{3^{m-1}p'}{3^{n}q'}}^{2}.\nonumber
\end{align}
\begin{remark} \label{maxk}
Using (\ref{eq:Mdef}) we infer that the maximal $k$ such that $n-(k-1)\in\M$
satisfies
\[
\frac{p}{3^{n-k+1}}\le\frac{1}{3}\text{ and }\frac{p}{3^{n-k}}>\frac{1}{3}.
\]
This is equivalent with
\begin{equation}
k-n=\left\lfloor \log_{3}\frac{1}{p}\right\rfloor =\left\lfloor \log_{3}\frac{q'}{p'}\right\rfloor .\label{eq:maxk}
\end{equation}
\end{remark}
\subsubsection{Quadraric variations along the Lebesgue partitions for the grids
$\protect\rbr{p/3^{n}}\protect\Z+\protect\rbr{r/3^{n}}$ when $r$
is irrational}

Formula (\ref{eq:k-step_rec}) simplifies when
\begin{equation}
\theta-\sum_{l=1}^{m}j_{l}\frac{3^{n}}{3^{l}}\frac{q'}{p'}\text{ for all }m=1,\ldots,k; j_{1},\ldots,j_{m-1}\in\cbr{0,1,2},j_{m}\in\cbr{1,2}, \text{ is not integer}.\label{eq:conddd}
\end{equation}
This holds for example when $r$ is irrational thus $\theta=r/p$ is an irrational number too.
In such a case for all $m=1,\ldots,k; j_{1},\ldots,j_{m-1}\in\cbr{0,1,2},j_{m}\in\cbr{1,2},$
we get
\[
\mathbf{1}_{(0,1)}\rbr{\left\{ \theta-\sum_{l=1}^{m}j_{l}\frac{3^{n}q'}{3^{l}p'}\right\} }=1,
\]
so the second sum in (\ref{eq:k-step_rec}) simplifies to
\begin{align}
 & \sum_{m=1}^{k}\frac{1}{3^{m-1}}\sum_{j_{1},...,j_{m-1}\in\cbr{0,1,2},j_{m}\in\cbr{1,2}}\mathbf{1}_{(0,1)}\rbr{\left\{ \theta-\sum_{l=1}^{m}j_{l}\frac{3^{n}q'}{3^{l}p'}\right\} }\rbr{\frac{3^{m-1}p'}{3^{n}q'}}^{2}\nonumber \\
 & =\sum_{m=1}^{k}\frac{1}{3^{m-1}}3^{m-1}2\rbr{\frac{3^{m-1}p'}{3^{n}q'}}^{2}=2\rbr{\frac{p'}{3^{n}q'}}^{2}\sum_{m=1}^{k}9^{m-1}\nonumber \\
 & =\frac{1}{4}\frac{9^{k}-1}{9^{n}}\rbr{\frac{p'}{q'}}^{2}\label{eq:sec_summand}
\end{align}
and (\ref{eq:k-step_rec}) simplifies to
\begin{align}
\sbr x_{1}^{c_{n},\theta c_{n}} & =\frac{1}{3^{k}}\sum_{i_{1},i_{2},\ldots,i_{k}\in\cbr{0,1,2}}\sbr x_{1}^{\frac{3^{k}p'}{3^{n}q'},\left\{ \theta-\sum_{v=1}^{k}i_{v}\frac{3^{n}q'}{3^{v}p'}\right\} \frac{3^{k}p'}{3^{n}q'}}+\frac{1}{4}\frac{9^{k}-1}{9^{n}}\rbr{\frac{p'}{q'}}^{2}.\label{eq:k-step_rec-1}
\end{align}

\begin{remark} \label{theeeta}
Notice that for $m=1,\ldots,k,\quad j_{1},\ldots,j_{m-1}\in\cbr{0,1,2},j_{m}\in\cbr{1,2},$
the number $\sum_{l=1}^{m}j_{l}\frac{3^{n}}{3^{l}}\frac{q'}{p'}$
may be represented as
\[
\frac{q'}{p'}3^{n}\frac{3^{m-1}j_{1}+\ldots3j_{m-1}+j_{m}}{3^{m}}=\frac{q'}{p'}3^{n-m}a,
\]
where $a\in\cbr{1,2,4,5,\ldots,3^{m}-2,3^{m}-1}$ or equivalently,
$a$ is an integer between $0$ and $3^{m}$, not divisible by $3$;
and opposite \textendash{} for any integer $a$ between $0$ and $3^{m}$,
not divisible by $3$, there exist numbers $j_{m-1}\in\cbr{0,1,2},j_{m}\in\cbr{1,2}$
such that
\[
a=3^{m-1}j_{1}+\ldots3j_{m-1}+j_{m}.
\]
Further, by (\ref{eq:maxk}), $m$ may be any positive integer no greater
than $k=n+\left\lfloor \log_{3}\frac{q'}{p'}\right\rfloor$ so
$n-m$ may be any integer no smaller than $-\left\lfloor \log_{3}\frac{q'}{p'}\right\rfloor $.
Since $n$ may be as large as we please (further we will tend with
$n$ to $+\ns$) thus $a$ may be any integer number not divisible
by $3$ and condition (\ref{eq:conddd}) is equivalent with the
condition that for any integer $b=n-m\ge-\left\lfloor \log_{3}\frac{q'}{p'}\right\rfloor $
and for any integer $a$ not divisible by $3$,
\begin{equation*}
\theta\neq\cbr{\frac{q'a}{p'}3^{b}}.
\end{equation*}
This is equivalent with the fact that $\theta$ is not a multiple of
\[
\frac{3^{-\left\lfloor \log_{3}\frac{q'}{p'}\right\rfloor }}{p'} = \frac{1}{3^{\left\lfloor \log_{3}\frac{q'}{p'}\right\rfloor }p'} .
\]
\end{remark}
\subsection{Limit of quadratic variations of $x$ along the Lebesgue partitions
for the grids $\protect\rbr{p/3^{n}}\protect\Z+\protect\rbr{r/3^{n}}$, 
$p\in\protect\Q_{+}$, $r\in\protect\R\setminus\protect\Q$ }

Now we are ready to prove the first main result of this paper. 

\begin{theorem} \label{main} Let $x:[0,1]\ra\R$ be the horizontal component of the
Peano curve, defined by (\ref{eq:px}). Let $p=p'/q'$ where $p'$
and $q'$ are positive integers not divisible by $3$ and $r$
be such that $r=\theta p$ where $\theta$ is not a multiple of any $3^{m}/p'$, $m\in \Z$.
Then for $t \in [0,1]$
\[
\lim_{n\ra+\ns}\sbr x_{t}^{3^{-n}p,3^{-n}r}=3^{\left\lfloor -\log_{3}p\right\rfloor }p\cbr{1-\frac{3}{4}\rbr{3^{\left\lfloor -\log_{3}p\right\rfloor }p}}t.
\]
\end{theorem}
{\bf Proof:} Let us fix $n\in\M$ and take the maximal $k=n+\left\lfloor \log_{3}\frac{q'}{p'}\right\rfloor $
such that $n-(k-1)\in\M$ so (\ref{eq:k-step_rec}) holds for any
real $\theta$. By Remark \ref{theeeta}, since $\theta$  is not a multiple of any $3^{m}/p'$, $m\in \Z$, \eqref{eq:k-step_rec-1} holds. Since the grids
\[
\frac{p}{3^{n}}\Z+\frac{r}{3^{n}}=\frac{p}{3^{n}}\rbr{\Z+\theta}\text{ and }\frac{p}{3^{n}}\rbr{\Z+\theta+q}
\]
are the same for any integer $q$, without loss of generality we may
assume that $\theta\in\left(0,1\right)$. There exists some $l=0,1,\ldots,3^{\max\cbr{k-n,0}}p'-1$
such that
\begin{equation}
\theta\in\left(\frac{l}{3^{\max\cbr{k-n,0}}p'},\frac{l+1}{3^{\max\cbr{k-n,0}}p'}\right).\label{eq:some_theta}
\end{equation}

To estimate $\sbr x_{1}^{c_{n},\theta c_{n}}$ we need to deal with
the sum:
\begin{equation}
\frac{1}{3^{k}}\sum_{i_{1},i_{2},\ldots,i_{k}\in\cbr{0,1,2}}\sbr x_{1}^{\frac{3^{k}p'}{3^{n}q'},\left\{ \theta-\sum_{v=1}^{k}i_{v}\frac{3^{n}q'}{3^{v}p'}\right\} \frac{3^{k}p'}{3^{n}q'}}.\label{eq:ssssum}
\end{equation}
Numbers $\sum_{v=1}^{k}i_{v}\frac{3^{n}q'}{3^{v}p'}=\frac{3^{n}q'}{3^{k}p'}\sum_{v=1}^{k}3^{k-v}i_{v}$,
$i_{1},i_{2},\ldots,i_{k}\in\cbr{0,1,2}$, appearing in \eqref{eq:ssssum},
range over $\frac{3^{n}q'}{3^{k}p'}M$, where $M\in\cbr{0,1,2,\ldots,3^{k}-1}$.

Let us notice that for $c\in(1/3,1]$ and $r\in[0,c)$ we have
\[
\sbr x_{1}^{c,r}=w\cdot c^{2}\text{ for }w=0,1,2
\]
iff
\begin{equation}
w\cdot c+r\le1\text{ and }(w+1)c+r>1\label{eq:wc}
\end{equation}
(notice that $3c>1$ thus $\sbr x_{1}^{c,r}$ can not be equal $3\rbr{\frac{3^{k}p'}{3^{n}q'}}^{2}$). For
example, if (\ref{eq:wc}) holds with $w=0$ this means that $x$
hits the net $c\Z+r$ on the time interval $[0,1]$ only at the points
$t\in[0,1]$ where $x(t)=r$ and $\sbr x_{1}^{c,r}=0$; if (\ref{eq:wc})
holds with $w=1$ this means that $x$ hits the net $c\Z+r$ only
at the points $t\in[0,1]$ such $x(t)=r$ or $x(t)=c+r$. Moreover,
if $x(t_{1})=r$ and $x(t_{2})=c+r$ then $t_{2}>t_{1}$. Thus $\sbr x_{1}^{c,r}=c^{2}$.

In particular, taking the maximal $k=n+\left\lfloor \log_{3}\frac{q'}{p'}\right\rfloor $,
$c=3^{k}p'/\rbr{3^{n}q'},$ $r=\left\{ \theta-3^{k}p'/M\rbr{3^{n}q'}\right\} c$,
we have
\[
\sbr x_{1}^{\frac{3^{k}p'}{3^{n}q'},\left\{ \theta-\frac{3^{n}q'}{3^{k}p'}M\right\} \frac{3^{k}p'}{3^{n}q'}}=w\rbr{\frac{3^{k}p'}{3^{n}q'}}^{2}\text{ for }w=0,1,2
\]
iff
\[
w\frac{3^{k}p'}{3^{n}q'}+\left\{ \theta-\frac{3^{n}q'}{3^{k}p'}M\right\} \frac{3^{k}p'}{3^{n}q'}\le1\text{ and }(w+1)\frac{3^{k}p'}{3^{n}q'}+\left\{ \theta-\frac{3^{n}q'}{3^{k}p'}M\right\} \frac{3^{k}p'}{3^{n}q'}>1.
\]

Let us consider two cases.

\textbf{1st case.}

\[
\frac{3^{k}p'}{3^{n}q'}>\frac{1}{2}.
\]
In this case we have
\[
\sbr x_{1}^{\frac{3^{k}p'}{3^{n}q'},\left\{ \theta-\frac{3^{n}q'}{3^{k}p'}M\right\} \frac{3^{k}p'}{3^{n}q'}}=\rbr{\frac{3^{k}p'}{3^{n}q'}}^{2}
\]
if and only if
\begin{equation}
\frac{3^{k}p'}{3^{n}q'}+\left\{ \theta-\frac{3^{n}q'}{3^{k}p'}M\right\} \frac{3^{k}p'}{3^{n}q'}\le1\label{eq:condition}
\end{equation}
and $\sbr x_{1}^{\frac{3^{k}p'}{3^{n}q'},\left\{ \theta-\frac{3^{n}q'}{3^{k}p'}M\right\} \frac{3^{k}p'}{3^{n}q'}}=0$
otherwise (notice that $2\frac{3^{k}p'}{3^{n}q'}>1$ thus $\sbr x_{1}^{\frac{3^{k}p'}{3^{n}q'},\left\{ \theta-\frac{3^{n}q'}{3^{k}p'}M\right\} \frac{3^{k}p'}{3^{n}q'}}$
can not be equal $2\rbr{\frac{3^{k}p'}{3^{n}q'}}^{2}$). Thus we need
to establish what is the number of $M\in\cbr{0,1,2,\ldots,3^{k}-1}$
such that (\ref{eq:condition}) holds. Let us notice that (\ref{eq:condition})
is equivalent with
\begin{equation}
\left\{ \theta-\frac{3^{n}q'}{3^{k}p'}M\right\} \le\frac{3^{n}q'}{3^{k}p'}-1.\label{eq:condition_eqv}
\end{equation}
Further, let us notice that
\[
0 < \frac{3^{n}q'}{3^{k}p'}-1=\frac{3^{\max\cbr{k-n,0}-(k-n)}q'-3^{\max\cbr{k-n,0}}p'}{3^{\max\cbr{k-n,0}}p'}<1
\]
and for non-negative integer $l$ satisfying (\ref{eq:some_theta}) there are exactly
\begin{equation} \label{numb_l}
3^{\max\cbr{k-n,0}-(k-n)}q'-3^{\max\cbr{k-n,0}}p'+1
\end{equation}
integers $L$ among $0,1,\ldots,3^{\max\cbr{k-n,0}}p'-1$ for
which
\[
\cbr{\frac{l}{3^{\max\cbr{k-n,0}}p'}-\frac{3^{n}q'}{3^{k}p'}L}= \cbr{\frac{l-3^{\max\cbr{k-n,0}-(k-n)}q'L}{3^{\max\cbr{k-n,0}}p'}}\le\frac{3^{n}q'}{3^{k}p'}-1,
\]
namely
\[
L=Ql,Q(l-1),\ldots, Q\rbr{l-\cbr{3^{\max\cbr{k-n,0}-(k-n)}q'-3^{\max\cbr{k-n,0}}p'}}\text{ (mod }3^{\max\cbr{k-n,0}}p'),
\]
where $Q$ is an inverse of $3^{\max\cbr{k-n,0}-(k-n)}q'$ $\text{ (mod }3^{\max\cbr{k-n,0}}p')$, that is $Q \in \cbr{1,2,\ldots, 3^{\max\cbr{k-n,0}}p'-1}$ and 
\[
3^{\max\cbr{k-n,0}-(k-n)}q'Q = 1 \text{ (mod }3^{\max\cbr{k-n,0}}p')
\] 
($a\text{ (mod }b$) is the nonnegative remainder from the division
of the integer $a$ by the integer $b$). 
Such a number $Q$ exists since $3^{\max\cbr{k-n,0}-(k-n)}q'$ and $3^{\max\cbr{k-n,0}}p'$ are relatively prime, see for example \cite[Theorem 57]{HW:1975}.

Analogously, the number of integers $M$ among $0,1,\ldots,3^{\max\cbr{k-n,0}}p'-1$
for which (\ref{eq:condition_eqv}), equivalently,
condition (\ref{eq:condition}) holds, is equal to
\begin{equation} \label{numbM}
N = 3^{\max\cbr{k-n,0}-(k-n)}q'-3^{\max\cbr{k-n,0}}p'.
\end{equation}
It is is smaller from \eqref{numb_l} by $1$ since $\theta > {l}/\rbr{3^{\max\cbr{k-n,0}}p'}$ and we can not take $$Q\rbr{l-\cbr{3^{\max\cbr{k-n,0}-(k-n)}q'-3^{\max\cbr{k-n,0}}p'}}\text{ (mod }3^{\max\cbr{k-n,0}}p').$$

Proportionally, there are no more that
\[
\rbr{3^{\max\cbr{k-n,0}-(k-n)}q'-3^{\max\cbr{k-n,0}}p'}\rbr{\frac{3^{k}}{3^{\max\cbr{k-n,0}}p'} +1} = \rbr{\frac{3^{n}q'}{3^{k}p'}-1}3^k+N
\]
and no less than
\[
\rbr{3^{\max\cbr{k-n,0}-(k-n)}q'-3^{\max\cbr{k-n,0}}p'}\rbr{\frac{3^{k}}{3^{\max\cbr{k-n,0}}p'} -1} = \rbr{\frac{3^{n}q'}{3^{k}p'}-1}3^k-N
\]
integers $M$ among $0,1,\ldots,3^{k}-1$ ($k$ is large) for
which
\[
\cbr{\theta-\frac{3^{n}q'}{3^{k}p'}M}\le\frac{3^{n}q'}{3^{k}p'}-1
\]
and
\begin{equation*}
\sum_{i_{1},i_{2},\ldots,i_{k}\in\cbr{0,1,2}}\sbr x_{1}^{\frac{3^{k}p'}{3^{n}q'},\left\{ \theta-\sum_{v=1}^{k}i_{v}\frac{3^{n}q'}{3^{v}p'}\right\} \frac{3^{k}p'}{3^{n}q'}} = \rbr{\rbr{\frac{3^{n}q'}{3^{k}p'}-1} 3^k+\nu}\rbr{\frac{3^{k}p'}{3^{n}q'}}^{2},
\end{equation*}
where $\nu \in [-N,N]$.
Recalling that $k=n+\left\lfloor \log_{3}\frac{q'}{p'}\right\rfloor $ we get
\begin{align*}
 & \lim_{n\ra+\ns}\frac{1}{3^{k}}\sum_{i_{1},i_{2},\ldots,i_{k}\in\cbr{0,1,2}}\sbr x_{1}^{\frac{3^{k}p'}{3^{n}q'},\left\{ \theta-\sum_{v=1}^{k}i_{v}\frac{3^{n}q'}{3^{v}p'}\right\} \frac{3^{k}p'}{3^{n}q'}}\\
 & =\lim_{n\ra+\ns}\frac{1}{3^{k}}\rbr{\frac{3^{k}p'}{3^{n}q'}}^{2}\rbr{\rbr{\frac{3^{n}q'}{3^{k}p'}-1} 3^k+\nu}\rbr{\frac{3^{k}p'}{3^{n}q'}}^{2}=3^{\left\lfloor \log_{3}\frac{q'}{p'}\right\rfloor }\rbr{\frac{p'}{q'}}-\rbr{3^{\left\lfloor \log_{3}\frac{q'}{p'}\right\rfloor }\frac{p'}{q'}}^{2}.
\end{align*}
Finally, by (\ref{eq:k-step_rec-1}) we have
\begin{align*}
\lim_{n\ra+\ns}\sbr x_{1}^{c_{n},\theta c_{n}} & =3^{\left\lfloor \log_{3}\frac{q'}{p'}\right\rfloor }\rbr{\frac{p'}{q'}}-\frac{3}{4}\rbr{3^{\left\lfloor \log_{3}\frac{q'}{p'}\right\rfloor }\frac{p'}{q'}}^{2}\\
 & =3^{\left\lfloor -\log_{3}p\right\rfloor }p-\frac{3}{4}\rbr{3^{\left\lfloor -\log_{3}p\right\rfloor }p}^{2}\\
 & =3^{\left\lfloor -\log_{3}p\right\rfloor }p\cbr{1-\frac{3}{4}3^{\left\lfloor -\log_{3}p\right\rfloor }p}.
\end{align*}

\textbf{2nd case.}
\[
\frac{3^{k}p'}{3^{n}q'}\le\frac{1}{2}.
\]

Similarly as in the first case we need to calculate for how many $M=0,1,2,\ldots,3^{k}-1$
one has
\begin{equation}
\frac{3^{k}p'}{3^{n}q'}+\left\{ \theta-\frac{3^{n}q'}{3^{k}p'}M\right\} \frac{3^{k}p'}{3^{n}q'}\le1\text{ and }2\frac{3^{k}p'}{3^{n}q'}+\left\{ \theta-\frac{3^{n}q'}{3^{k}p'}M\right\} \frac{3^{k}p'}{3^{n}q'}>1,\label{eq:first_cond}
\end{equation}
or, equivalently, for how many $M=0,1,2,\ldots,3^{k}-1$ one has
\[
\left\{ \theta-\frac{3^{n}q'}{3^{k}p'}M\right\} \le\frac{3^{n}q'}{3^{k}p'}-1\text{ and }\left\{ \theta-\frac{3^{n}q'}{3^{k}p'}M\right\} >\frac{3^{n}q'}{3^{k}p'}-2
\]
and secondly, for how many $M=0,1,2,\ldots,3^{k}-1$ one has
\begin{equation}
2\frac{3^{k}p'}{3^{n}q'}+\left\{ \theta-\frac{3^{n}q'}{3^{k}p'}M\right\} \frac{3^{k}p'}{3^{n}q'}\le1.\label{eq:second_cond}
\end{equation}
Reasoning similarly as in the first case we get that there are 
$\rbr{\frac{3^{n}q'}{3^{k}p'}-2}3^{k}+\nu_1$ integers $M$ among numbers
$0,1,\ldots,3^{k}-1$ ($k$ large) for which (\ref{eq:second_cond})
holds and $\rbr{1-\rbr{\frac{3^{n}q'}{3^{k}p'}-2}}3^{k}+\nu_2$ integers
$M$ among numbers $0,1,\ldots,3^{k}-1$ ($k$ large) for which (\ref{eq:first_cond})
holds; $\nu_1, \nu_1 \in \sbr{-N,N}$.

As a result we get also that
\begin{align*}
 & \lim_{n\ra+\ns}\frac{1}{3^{k}}\sum_{i_{1},i_{2},\ldots,i_{k}\in\cbr{0,1,2}}\sbr x_{1}^{\frac{3^{k}p'}{3^{n}q'},\left\{ \theta-\sum_{v=1}^{k}i_{v}\frac{3^{n}q'}{3^{v}p'}\right\} \frac{3^{k}p'}{3^{n}q'}}\\
 & =\lim_{n\ra+\ns}\sbr{\frac{2}{3^{k}}  \rbr{\frac{3^{k}p'}{3^{n}q'}}^{2} \rbr{ \rbr{\frac{3^{n}q'}{3^{k}p'}-2}3^{k} +\nu_1}+\frac{1}{3^{k}}\rbr{\frac{3^{k}p'}{3^{n}q'}}^{2} \rbr{ \rbr{3-\frac{3^{n}q'}{3^{k}p'}}3^{k} + \nu_2}  }\\
 & =3^{\left\lfloor \log_{3}\frac{q'}{p'}\right\rfloor }\rbr{\frac{p'}{q'}}-\rbr{3^{\left\lfloor \log_{3}\frac{q'}{p'}\right\rfloor }\frac{p'}{q'}}^{2}
\end{align*}
and by (\ref{eq:k-step_rec-1}) we have
\begin{align*}
\lim_{n\ra+\ns}\sbr x_{1}^{c_{n},\theta c_{n}} & =3^{\left\lfloor \log_{3}\frac{q'}{p'}\right\rfloor }\rbr{\frac{p'}{q'}}-\frac{3}{4}\rbr{3^{\left\lfloor \log_{3}\frac{q'}{p'}\right\rfloor }\frac{p'}{q'}}^{2}\\
 & =3^{\left\lfloor -\log_{3}p\right\rfloor }p\cbr{1-\frac{3}{4}\rbr{3^{\left\lfloor -\log_{3}p\right\rfloor }p}}.
\end{align*}
Let us notice that the obtained limit is the same as in the first case.

Now we are going to prove convergence of $\sbr x_{t}^{c_{n},\theta c_{n}} $ for any $t \in [0,1]$.
Using the relations \eqref{eq:s1}-\eqref{eq:s9}, we infer that
\[
\lim_{n\ra+\ns}\sbr x_{i/9}^{3^{-n}p,3^{-n}r}=3^{\left\lfloor -\log_{3}p\right\rfloor }p\cbr{1-\frac{3}{4}\rbr{3^{\left\lfloor -\log_{3}p\right\rfloor }p}}\frac{i}{9}, \quad i=1,2,\ldots,9.
\]
Similarly, using multiple times  \eqref{eq:s1}-\eqref{eq:s9}, we get
\[
\lim_{n\ra+\ns}\sbr x_{i/9^m}^{3^{-n}p,3^{-n}r}=3^{\left\lfloor -\log_{3}p\right\rfloor }p\cbr{1-\frac{3}{4}\rbr{3^{\left\lfloor -\log_{3}p\right\rfloor }p}}\frac{i}{9^m}, \quad i=1,2,\ldots,9^m, m=1,2,\ldots.
\]
From the fact that $[0,1] \ni t \mapsto \liminf_{n\ra+\ns}\sbr x_{t}^{3^{-n}p,3^{-n}r}$ and $[0,1] \ni t \mapsto \limsup_{n\ra+\ns}\sbr x_{t}^{3^{-n}p,3^{-n}r}$ are not decreasing in $t$, and from two last limits we get that
\[
\lim_{n\ra+\ns}\sbr x_{t}^{3^{-n}p,3^{-n}r}=3^{\left\lfloor -\log_{3}p\right\rfloor }p\cbr{1-\frac{3}{4}\rbr{3^{\left\lfloor -\log_{3}p\right\rfloor }p}}t, \quad t\in [0,1].
\]
\hfill $\square$

\subsection{Truncated variation and numbers of interval crossings of the $x$--component of the Peano curve} \label{sect_tv}

Let us introduce the following functional called \emph{truncated variation}.
The truncated variation of $x$ with the truncation parameter $c \ge 0$ on the time interval $[s,t]$, $0 \le s < t \le 1$ is defined as:
\begin{equation} \label{tvc_def}
\TTV x{\left[s,t\right]}{c}:=\sup_{\pi \in \Pi(s,t)}\sum_{[u,v] \in \pi} \max \rbr{\left| x_{v}-x_{u} \right| -c,0},
\end{equation}
where the supremum is taken over all finite partitions $\pi$ of the interval $[s,t]$, that is finite sets of no overlapping (with disjoint interiors) subintervals $[u,v]$ of $[s,t]$ such that $\bigcup_{[u,v] \in \pi}[u,v] = [s,t]$. The family of all such partitions is denoted by $\Pi(s,t)$.
The truncated variation is finite for any c\`adl\`ag or even regulated (i.e. possessing right- and left- limits) path $x:[0, 1] \ra \R$ whenever $c>0$.

Using \cite[last formula in Sect. 3.3]{BDL:2024}, Theorem \ref{main} (for example in the case when $r$ is irrational) and denoting
\[
C_p := 3^{\left\lfloor -\log_{3}p\right\rfloor }p\cbr{1-\frac{3}{4}\rbr{3^{\left\lfloor -\log_{3}p\right\rfloor }p}}
\]
we get that for the $x$--component of the Peano curve we have the following convergence
\begin{equation} \label{tv}
\lim_{n \ra +\ns} 3^{-n}p \TTV x{\left[0,t\right]}{3^{-n}p} = C_p t, \quad t\in [0,1].
\end{equation}

Now, we introduce \emph{numbers of interval (up- and down-) crossings} by continuous $x:[0,1] \ra \R$ on the time interval $[s, t]$.
Let $z\in \R, c > 0$, define $\sigma_{0}^{c}=s$ and for $n \in \N_0$
\[
\tau_{n}^{c}=\inf\left\{ u \in \sbr{ \sigma_{n}^{c}, t}: x_u\ge z+c/2\right\} ,
\]
\[
\sigma_{n+1}^{c}=\inf\left\{ u \in \sbr{ \tau_{n}^{c}, t} : x_u < z-c/2\right\},
\]
where we apply conventions: $\inf \emptyset = +\ns$, $\sbr{+\ns, t} = \emptyset$.
\begin{dfn}\label{defd}
The \emph{number of downcrossings, upcrossings and crossings  by $x$ the interval $[z-c/2, z+c/2]$} on the time interval $[s,t]$ one defines respectively as
\[
\Dcrossemph{z,c}{x}{[s,t]} :=\max\left\{ n:\sigma_{n}^{c}\leq t\right\}, \quad\Ucrossemph{z,c}{x}{[s,t]} := \Dcrossemph{-z,c}{-x}{[s,t]}.
\]
and
\[
\crossemph{z,c}{x}{[s,t]}:= \Ucrossemph{z,c}{x}{[s,t]} + \Dcrossemph{z,c}{x}{[s,t]}.
\]
\end{dfn}

Using \eqref{tv}, \cite[Theorem 2.2]{BDL:2024} and \cite[Remark 3.9]{BDL:2024} we get that for any continuous $g: \R \ra \R$,
\begin{equation} \label{conv_tv_zeta}
3^{-n}p \int_{\R} \cross{z,3^{-n}p}{x}{[0,\cdot]} g(z)\dd z\rightarrow C_p \int_{(0,\cdot]}g(x_{t})\dd t.
\end{equation}

\section{Local times of the $x$--component of the Peano curve}

In this section we will prove existence of local time of $x$ with respect to the Lebesgue measure, that is, existence
of a function
\[
L:[0,1]\times\R\ra[0,+\ns)
\]
such that for any $t\in[0,1]$, $L_{t}:=L(t,\cdot):\R\ra[0,+\ns)$
is Borel-measurable and for any Borel-measurable and bounded function
$g:\R\ra\R$ one has
\begin{equation}
\int_{0}^{t}g\rbr{x_{s}}\dd s=\int_{\R}g(z)L_{t}^{z}\dd z.\label{eq:occ_times_formula}
\end{equation}
(In this section, for typographical reasons, we will rather use the notation $x_s$ instead of $x(s)$.)

\begin{theorem}
There exists  local time of the horizontal component $x$ of the Peano curve. Moreover, this local time attains values no greater than $1$.
\end{theorem}
{\bf Proof:} To prove (\ref{eq:occ_times_formula}) it is sufficient to prove that
for any interval $[a,b]\subset\R$, $a<b$, one has
\[
\int_{0}^{t}{\bf 1}_{[a,b]}\rbr{x_{s}}\dd s=\int_{a}^{b}L_{t}^{z}\dd z.
\]
We will construct $L$ recursively, and from the construction we will infer that $L$ attains values no greater than $1$.

Let us fix $n=1,2,\ldots$, such that $b-a>3^{-n}$, and approximate
$a$ and $b$ by numbers of the form $p_{a}/3^{n}$, $p_{b}/3^{n}$ resp.,
where $p_{a}$ and $p_{b}$ are integers satisfying
\[
\frac{p_{a}}{3^{n}}\le a<\frac{p_{a}+1}{3^{n}}\le\frac{p_{b}}{3^{n}}\le b < \frac{p_{b}+1}{3^{n}}.
\]

Using this approximation we have that it is sufficient to prove that
there exists a function
\[
L:[0,1]\times\R\ra[0,1]
\]
such that for any $t\in[0,1]$, $L_{t}:=L(t,\cdot):\R\ra[0,1]$
is Borel-measurable and for any numbers of the form
\[
a=\frac{p_{a}}{3^{n}},\quad b=\frac{p_{b}}{3^{n}} \ge \frac{p_{a}+1}{3^{n}},\quad n\in\N,\quad p_{a},p_{b}\in\Z,
\]
one has
\begin{equation} \label{ind_step}
\int_{0}^{t}{\bf 1}_{[a,b]}\rbr{x_{s}}\dd s=\int_{a}^{b}L_{t}^{z}\dd z.
\end{equation}

First, we will define $L_{t}$ for $t$ of the form $k\cdot9^{-N}$,
$N=1,2,\ldots$, $k=1,2,\ldots,9^{N}$.
{For any $t \in [0, 1]$ and $z \notin (0,1)$ we define $L^z_t=0$.}
Next we start with $N=1$, that is with $t$ of the form $u\cdot9^{-1}$,
$N=1,2,\ldots$, $u=1,2,\ldots,9$. Let us consider intervals $[a,b]$
of the form
\[
[a,b]=\sbr{\frac{p_{a}}{3},\frac{p_{a}+1}{3}},\quad p_{a}=0,1,2.
\]
Let us notice that
\begin{equation}
\int_{0}^{1/9}{\bf 1}_{[a,b]}\rbr{x_{s}}\dd s=\begin{cases}
\frac{1}{9} & \text{if }[a,b]=\sbr{\frac{0}{3},\frac{1}{3}},\\
0 & \text{if }[a,b]=\sbr{\frac{1}{3},\frac{2}{3}},\sbr{\frac{2}{3},\frac{3}{3}}
\end{cases}=\int_{a}^{b}L_{1/9}^{z}\dd z\label{eq:local_first}
\end{equation}
where
\[
L_{1/9}^{z}=\begin{cases}
\frac{1}{3} & \text{if }z\in\rbr{0,\frac{1}{3}},\\
0 & \text{if }z\in\rbr{\frac{1}{3},\frac{3}{3}}.
\end{cases}
\]
Equality (\ref{eq:local_first}) holds since for $s\in\sbr{\frac{0}{9},\frac{1}{9}}$
one has $x_{s}\in\sbr{\frac{0}{3},\frac{1}{3}}$.

Further we have
\[
\int_{0}^{2/9}{\bf 1}_{[a,b]}\rbr{x_{s}}\dd s=\begin{cases}
\frac{2}{9} & \text{if }[a,b]=\sbr{\frac{0}{3},\frac{1}{3}},\\
0 & \text{if }[a,b]=\sbr{\frac{1}{3},\frac{2}{3}},\sbr{\frac{2}{3},\frac{3}{3}}
\end{cases}=\int_{a}^{b}L_{2/9}^{z}\dd z,
\]
\[
\int_{0}^{3/9}{\bf 1}_{[a,b]}\rbr{x_{s}}\dd s=\begin{cases}
\frac{3}{9} & \text{if }[a,b]=\sbr{\frac{0}{3},\frac{1}{3}},\\
0 & \text{if }[a,b]=\sbr{\frac{1}{3},\frac{2}{3}},\sbr{\frac{2}{3},\frac{3}{3}}
\end{cases}=\int_{a}^{b}L_{3/9}^{z}\dd z,
\]
where
\[
L_{2/9}^{z}=\begin{cases}
\frac{2}{3} & \text{if }x\in\rbr{0,\frac{1}{3}},\\
0 & \text{if }z\in\rbr{\frac{1}{3},\frac{3}{3}},
\end{cases}\quad L_{3/9}^{z}=\begin{cases}
\frac{3}{3} & \text{if }z\in\rbr{0,\frac{1}{3}},\\
0 & \text{if }z\in\rbr{\frac{1}{3},\frac{3}{3}}.
\end{cases}
\]
Similarly, defining
\begin{align*}
L_{4/9}^{z} & =\begin{cases}
\frac{3}{3} & \text{if }z\in\rbr{0,\frac{1}{3}},\\
\frac{1}{3} & \text{if }z\in\rbr{\frac{1}{3},\frac{2}{3}},\\
0 & \text{if }z\in\rbr{\frac{2}{3},\frac{3}{3}};
\end{cases}\quad L_{5/9}^{z}=\begin{cases}
\frac{3}{3} & \text{if }z\in\rbr{0,\frac{1}{3}},\\
\frac{2}{3} & \text{if }z\in\rbr{\frac{1}{3},\frac{2}{3}},\\
0 & \text{if }z\in\rbr{\frac{2}{3},\frac{3}{3}};
\end{cases}\quad L_{6/9}^{z}=\begin{cases}
\frac{3}{3} & \text{if }z\in\rbr{0,\frac{1}{3}}\cup\rbr{\frac{1}{3},\frac{2}{3}},\\
0 & \text{if }z\in\rbr{\frac{2}{3},\frac{3}{3}};
\end{cases}
\end{align*}
\[
L_{7/9}^{z}=\begin{cases}
\frac{3}{3} & \text{if }z\in\rbr{0,\frac{1}{3}}\cup\rbr{\frac{1}{3},\frac{2}{3}},\\
\frac{1}{3} & \text{if }z\in\rbr{\frac{2}{3},\frac{3}{3}};
\end{cases}\quad L_{8/9}^{z}=\begin{cases}
\frac{3}{3} & \text{if }z\in\rbr{0,\frac{1}{3}}\cup\rbr{\frac{1}{3},\frac{2}{3}},\\
\frac{2}{3} & \text{if }z\in\rbr{\frac{2}{3},\frac{3}{3}};
\end{cases}\quad L_{9/9}^{z}=1\text{ if }z\in(0,1)
\]
we get \eqref{ind_step} for $t=k/9$, $k=4,5,\ldots,9$, and $[a,b]=\sbr{{p_{a}}/{3},\rbr{p_{a}+1}/{3}}$, $p_{a}=0,1,2$. Recalling that for any $t \in [0, 1]$ and $z \notin (0,1)$ we have already defined $L^z_t=0$, we get \eqref{ind_step} for for $t=k/9$, $k=1,2,\ldots,9$ and  $[a,b]=\sbr{{p_{a}}/{3},p_{b}/{3}}$, $p_{a}, p_b \in \Z$, $p_a< p_b$.

Having defined $L_{t}$ for $t=k/9^{N}$, $k=1,2,\ldots, 9^{N}$,
we can define $L_{t}$ for $t$ from a finer grid, namely for $t=k/9^{N+1}$, $k=1,2,\ldots, 9^{N+1}$,
and check that for such $t$ equalities \eqref{ind_step}
hold for any $[a,b]$ of the form $a={p_{a}}/{3^{n}}, b={p_{b}}/{3^{n}}=\rbr{p_{a}+1}/{3^{n}},\ n=N+1$, $p_a \in \Z$ 
(hence for any $[a,b]$ of the form $a={p_{a}}/{3^{n}},\ b={p_{b}}/{3^{n}}\ge \rbr{p_{a}+1}/{3^{n}},\ n=0,1,2,\ldots,N+1$, $p_a, p_b \in \Z$).

Indeed, having defined $L_{t}$ for $t=k/9^{N}$, $k=1,2,\ldots, 9^{N}$ and assuming that \eqref{ind_step} holds for such $t$ and any $[a,b]$ of the form $a={p_{a}}/{3^{n}},\ b={p_{b}}/{3^{n}}\ge \rbr{p_{a}+1}/{3^{n}},\ n=0,1,2,\ldots,N$, $p_a, p_b \in \Z$,  let us define
\begin{align} \label{locallll}
L_{t}^{z}=\begin{cases}
\frac{1}{3} L_{9t}^{3z} & \text{if } t=\frac{1}{9^{N+1}},\frac{2}{9^{N+1}},\ldots,\frac{9^N}{9^{N+1}},\\
 L_{1/9}^{z} + \frac{1}{3} L_{9t-1}^{1 - 3z}  & \text{if } t=\frac{9^N+1}{9^{N+1}},\frac{9^N+2}{9^{N+1}},\ldots,\frac{2\cdot9^N}{9^{N+1}},\\
 L_{2/9}^{z} + \frac{1}{3} L_{9t-2}^{3z}  & \text{if } t=\frac{2 \cdot 9^N+1}{9^{N+1}},\frac{2 \cdot  9^N+2}{9^{N+1}},\ldots,\frac{3\cdot9^N}{9^{N+1}}, \\
L_{3/9}^{z} + \frac{1}{3} L_{9t-3}^{3z-1}  & \text{if } t=\frac{3 \cdot 9^N+1}{9^{N+1}},\frac{3 \cdot  9^N+2}{9^{N+1}},\ldots,\frac{4\cdot9^N}{9^{N+1}}, \\
L_{4/9}^{z} + \frac{1}{3} L_{9t-4}^{2-3z}  & \text{if } t=\frac{4 \cdot 9^N+1}{9^{N+1}},\frac{4 \cdot  9^N+2}{9^{N+1}},\ldots,\frac{5\cdot9^N}{9^{N+1}}, \\
L_{5/9}^{z} + \frac{1}{3} L_{9t-5}^{3z-1}  & \text{if } t=\frac{5 \cdot 9^N+1}{9^{N+1}},\frac{5 \cdot  9^N+2}{9^{N+1}},\ldots,\frac{6\cdot9^N}{9^{N+1}}, \\
 L_{6/9}^{z} + \frac{1}{3} L_{9t-6}^{3z-2}  & \text{if } t=\frac{6 \cdot 9^N+1}{9^{N+1}},\frac{6 \cdot  9^N+2}{9^{N+1}},\ldots,\frac{7\cdot9^N}{9^{N+1}}, \\
L_{7/9}^{z} + \frac{1}{3} L_{9t-7}^{3-3z}  & \text{if } t=\frac{7 \cdot 9^N+1}{9^{N+1}},\frac{7 \cdot  9^N+2}{9^{N+1}},\ldots,\frac{8\cdot9^N}{9^{N+1}}, \\
L_{8/9}^{z} + \frac{1}{3} L_{9t-8}^{3z-2}  & \text{if } t=\frac{8 \cdot 9^N+1}{9^{N+1}},\frac{8 \cdot  9^N+2}{9^{N+1}},\ldots,\frac{9\cdot9^N}{9^{N+1}}. \\
 \end{cases}
\end{align}

Now, for $t=k/9^{N+1}$, $k=1,2,\ldots, 9^N$, using \eqref{eq:s1}, substituting $u = 9s$, using the induction assumption and then substituting $w = 3z$ we obtain that for any $[a,b]$ of the form $a={p_{a}}/{3^{n}}, b={p_{b}}/{3^{n}} \ge \rbr{p_{a}+1}/{3^{n}},\ n=N+1$, $p_a, p_b \in \Z$,
\begin{align*}
\int_{0}^{t}{\bf 1}_{[a,b]}\rbr{x_{s}}\dd s & =\int_{0}^{t}{\bf 1}_{[a,b]}\rbr{\frac{1}{3}x_{9s}}\dd s = \frac{1}{9} \int_{0}^{9t}{\bf 1}_{[3a,3b]}\rbr{x_{u}}\dd u \\
& = \frac{1}{9} \int_{3a}^{3b} L_{9t}^{w}\dd w = \frac{1}{3} \int_{a}^{b} L_{9t}^{3z} \dd z,
\end{align*}
hence, setting $L_{t}^{z} = \frac{1}{3} L_{9t}^{3z}$ for $t=k/9^{N+1}$, $k=1,2,\ldots, 9^N$, we get \eqref{ind_step}.

Also, taking for example $t=k/9^{N+1}$, $k=4\cdot 9^{N} +1, 4\cdot 9^{N} + 2,\ldots, 5\cdot 9^{N}$, and assuming that \eqref{ind_step} holds for $u=m/9^{N+1}$, $m=1,2,\ldots, 4\cdot 9^N$, and any $[a,b]$ of the form $a={p_{a}}/{3^{n}}, b={p_{b}}/{3^{n}}\ge \rbr{p_{a}+1}/{3^{n}},\ n=N+1$, $p_a, p_b \in \Z$, we obtain
\begin{align}
\int_{0}^{t}{\bf 1}_{[a,b]}\rbr{x_{s}}\dd s & =\int_{0}^{4/9}{\bf 1}_{[a,b]}\rbr{x_{s}}\dd s + \int_{4/9}^{t}{\bf 1}_{[a,b]}\rbr{x_{s}}\dd s \nonumber  \\
& = \int_{a}^{b} L_{4/9}^z \dd z + \int_{4/9}^{t}{\bf 1}_{[a,b]}\rbr{x_{s}}\dd s. \label{ind1}
\end{align}
Further, using \eqref{eq:s5} and substituting $u = 5-9s$ and then using \eqref{eq:s10} and substituting $v = 1-u$  we get
\begin{align}
\int_{4/9}^{t}{\bf 1}_{[a,b]}\rbr{x_{s}}\dd s & =\int_{4/9}^{t}{\bf 1}_{[a,b]}\rbr{\frac{1}{3} + \frac{1}{3}x_{5-9s}}\dd s = 
- \frac{1}{9} \int_{1}^{5-9t}{\bf 1}_{[a,b]}\rbr{\frac{1}{3} +\frac{1}{3}x_{u}}\dd u \nonumber \\
& =  \frac{1}{9} \int_{5-9t}^{1}{\bf 1}_{[3a-1,3b-1]}\rbr{x_{u}}\dd u  =  \frac{1}{9} \int_{5-9t}^{1}{\bf 1}_{[3a-1,3b-1]}\rbr{1-x_{1-u}}\dd u \nonumber \\
& =  \frac{1}{9} \int_{0}^{9t-4}{\bf 1}_{[2-3b,2-3a]}\rbr{x_{v}}\dd v =  \frac{1}{9} \int_{2-3b}^{2-3a} L_{9t-4}^w \dd w  = \frac{1}{3} \int_{a}^{b} L_{9t-4}^{2-3z}\dd z . \label{ind2}
\end{align}
We also used the induction assumption and the substitution $2-3z  = w$.
From \eqref{ind1} and \eqref{ind2} we get for  $t=k/9^{N+1}$, $k=4\cdot 9^{N} +1, 4\cdot 9^{N} + 2,\ldots, 5\cdot 9^{N}$, and any $[a,b]$ of the form $a={p_{a}}/{3^{n}}, b={p_{b}}/{3^{n}} \ge \rbr{p_{a}+1}/{3^{n}},\ n=N+1$, $p_a, p_b \in \Z$,
\[
\int_{0}^{t}{\bf 1}_{[a,b]}\rbr{x_{s}}\dd s = \int_{a}^{b} L_{4/9}^z + \frac{1}{3} L_{9t-4}^{2-3z} \ \dd z,
\]
hence, setting $L_{t}^{z} = L_{4/9}^z + \frac{1}{3} L_{9t-4}^{2-3z}$ for $t=k/9^{N+1}$, $k=4\cdot 9^N + 1,4\cdot 9^N + 2,\ldots, 5 \cdot 9^{N}$, we get \eqref{ind_step}.

In a similar way, using \eqref{eq:s2}-\eqref{eq:s6}, \eqref{eq:s8}-\eqref{eq:s9} and \eqref{eq:s10}, we can prove that \eqref{locallll} defines recursively $L_t$ for all remaining $t$ of the form $t=k/9^{N+1}$, $k=1,2,\ldots, 9^{N+1}$,
such that for such $t$ equalities \eqref{ind_step} 
hold for any $[a,b]$ of the form $a={p_{a}}/{3^{n}},\ b={p_{b}}/{3^{n}}\ge \rbr{p_{a}+1}/{3^{n}},\ n=0,1,2,\ldots,N+1$, $p_a, p_b \in \Z$.
As a result, proceeding with $N$ to $+\ns$ we obtain $L_t$ defined for any $t$ of the form $k/9^N$, $k =1,2,\ldots, 9^N$, $N \in \N$,  and such that \eqref{ind_step} 
holds for any $[a,b]$ of the form $a={p_{a}}/{3^{n}},\ b={p_{b}}/{3^{n}}\ge \rbr{p_{a}+1}/{3^{n}},\ n\in \N$, $p_a, p_b \in \Z$.

Moreover, using values of $L_t^z$ for $t=1/9,2/9,\ldots, 8/9,1$, and the fact that $L_t^z = 0$ for any $t\in [0,1]$ and $z \notin (0,1)$ we infer inductively that $L_t^z$ defined by \eqref{locallll} attains values from the interval $[0,1]$.

\hfill $\square$

\subsection{Local time of the $x$-component of the Peano curve expressed as the weak limit of normalized numbers of interval crossings}
Recall the number $\cross{z,c}{x}{[0,t]}$ of crossing by $x$ the interval $[z-c/2, z+c/2]$ on the time interval $[0,t]$. In this section we will prove that the local time $L$ of the $x$-component of the Peano curve may be expressed as the weak limit of normalized numbers of interval crossings by $x$. This means that there exists some  $\varphi:(0,1) \ra [0,+\ns)$ such that $\lim_{c \ra 0+} \varphi(c) = 0$ and for any continuous $g: \R \ra \R$ 
\begin{equation} \label{loc_phi}
\lim_{c\ra 0+}  \int_{\R} \varphi(c) \cross{z,c}{x}{[0,t]} g(z)\dd z = \int_{\R} L_t^z  g(z)\dd z.
\end{equation}

Recall that \eqref{eq:occ_times_formula} states that
\[
\int_{\R}g(z) L_t^z  \dd z = \int_{0}^{t}g\rbr{x_{s}}\dd s.
\]

Now we will find  $\varphi:(0,1) \ra [0,+\ns)$ such that $\lim_{c \ra 0+} \varphi(c) = 0$ and 
\begin{equation} \label{loc_phiii}
\lim_{c\ra 0+}  \int_{\R} \varphi(c) \cross{z,c}{x}{[0,t]} g(z)\dd z = \int_{0}^{t}g\rbr{x_{s}}\dd s.
\end{equation}
Indeed, using \eqref{conv_tv_zeta} for $p>0$ we get
\begin{equation*}
 C_p^{-1} \int_{\R} 3^{-n}p \cdot \cross{z,3^{-n}p}{x}{[0,t]} g(z)\dd z\rightarrow \int_{0}^tg(x_{s})\dd s,
\end{equation*}
where $
C_p := 3^{\left\lfloor -\log_{3}p\right\rfloor }p\cbr{1-\frac{3}{4}\rbr{3^{\left\lfloor -\log_{3}p\right\rfloor }p}}
$.
Noticing that $C_p = C_{3^kp}$ for any integer $k$, we get
\begin{equation*}
 C_{3^{-n}p}^{-1} \int_{\R} 3^{-n}p \cdot \cross{z,3^{-n}p}{x}{[0,t]} g(z)\dd z\rightarrow \int_{0}^tg(x_{s})\dd s.
\end{equation*}
Thus \eqref{loc_phiii} is satisfied with $$\varphi(c) :=  C_{c}^{-1} \cdot c = \frac{1}{3^{\left\lfloor -\log_{3}c\right\rfloor }\cbr{1-\frac{3}{4}\rbr{3^{\left\lfloor -\log_{3}c\right\rfloor }c}}} .$$
The graph of the function $\varphi(c)$ (blue) together with graphs of $3c$ and $4c$ is presented on the figure below.
\begin{center}
\includegraphics[scale=0.65]{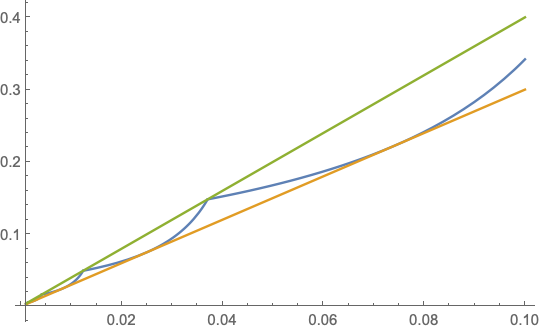}
\captionof{figure}{The graph of the function $\varphi(c)$ (blue) together with graphs of $3c$ and $4c$.}
\end{center}
Using \eqref{eq:occ_times_formula} and \eqref{loc_phiii} we finally get \eqref{loc_phi}.
\section*{Acknowledgments} The work of PLZ, DH and FJM was supported by a University Staff Doctoral Programme: Building Capacity in Applied Mathematics (USDP-BCAM) grant, ID 32, under the Newton's Operational Development Assistance Fund. The grant is funded by the UK Department for Business, Energy and Industrial Strategy and the Department of Higher Education and Training, South Africa and delivered by the British Council. For further information, please visit https://www.ukri.org/what-we-do/browse-our-areas-of-investment-and-support/newton-fund/.

The research of RM{\L} was supported by grant no. 2022/47/B/ST1/02114 \emph{Non-random equivalent characterizations of sample boundedness} of National Science Centre, Poland. 


\bibliographystyle{imsart-number}

\end{document}